\documentclass[papera4,12pt]{article}

\usepackage[english]{babel}
\usepackage[letterpaper,top=2cm,bottom=2cm,left=3cm,right=3cm,marginparwidth=1.75cm]{geometry}
\usepackage[utf8]{inputenc}

\usepackage{url}
\usepackage{amsmath}
\usepackage{graphicx}
\usepackage[english]{babel}
\usepackage[colorlinks=true, allcolors=red]{hyperref }
\usepackage{amsthm}
\usepackage{amsfonts}
\usepackage{hyperref}
\newtheorem{theorem}{Theorem}[section]
\newtheorem{lemma}{Lemma}[section]

\numberwithin{equation}{section}
\begin{document}

\begin{center}
{\bf {\Large On the Diophantine Equations

$J_n +J_m =L_k$ 

and

$L_n +L_m =J_k$ 

}}
\vspace{3mm}
\author[\hspace{0.7cm}\centerline{O. Salah$^*$ $ ^ 1$, A. Elsonbaty$^*$ $ ^ 2$, M. Anwar $^*$ $ ^ 3$ and M. A. Seoud $^*$ $ ^ 4$ }

$^*$Department of Mathematics, faculty of science, Ain Shams university\\
\indent \,\,\, e-mail: osamasalah@sci.asu.edu.eg$^*$ $ ^ 1$

\indent \,\,\, Ahmadelsonbaty@sci.asu.edu.eg$^*$ $ ^ 2$,\\ 
\indent \,\,\, mohamedanwar@sci.asu.edu.eg$^*$ $ ^ 3$,\\
\indent \,\,\, m.a.seoud@hotmail.com$^*$ $ ^ 4$.



\end{center}

\begin{abstract}
This study discovers Lucas numbers that are sum of two Jacobsthal numbers and every  Jacobsthal numbers that are sum of two Lucas numbers. Generally, we find all non-negative integer solutions $(n, m, k)$ of the two Diophantine equations $L_n +L_m =J_k$ and $J_n +J_m =L_K,$ where  $\left\lbrace L_{k}\right\rbrace_{k\geq0}$ and $\left\lbrace J_{n}\right\rbrace_{n\geq0}$ are the sequences of Lucas and Jacobsthal numbers, respectively. Our main findings are suppoorted by a modification of the Baker's theorem for linear forms in logarithms and  Dujella and Peth\H{o}'s reduction method.
\end{abstract}
\begin{flushleft}
\end{flushleft}
\begin{flushleft}
\end{flushleft}
\section{Introduction}
The Lucas sequence characterised by $L_0=2$, $L_1=1$ and 
\begin{equation*}
L_{n+2}=L_{n+1}+L_{n} \;\text{for all}\; n\geq0 . 
\end{equation*}
The history and applications of the Lucas sequence can be found in 
 \cite{Koshy2019}.
The Jacobsthal sequence characterised by $J_0=0$, $J_1=1$ and
\begin{equation*}
J_{n+2}=J_{n+1}+2 J_{n} \;\text{for all}\; n\geq0 . 
\end{equation*}
Key characteristics and abstract ideas of Jacobsthal numbers can be found in \cite{Horadam1996, Horadam1997, Horadam1988,Djordjevic2000, Yilmaz2009, Aydin2018}.

Finding all the terms in some integer sequences—which are the sums, products, and differences between two terms in another integer sequence have been considered in the literature by many researchers. For instance, Bravo and Luca found a solution to express powers of two as the sum of two Fibonacci numbers or as the sum of two Lucas numbers in \cite{Bravo2014,Bravo2016}, respectively. In 2020, Siar and Keskin\cite{Siar2017} solved the Diophantine equation $F_n-F_m=2^k$. Additionally, the sums of two Padovan numbers, or Fibonacci numbers, were all found in \cite{Garcia2022}. In 2021, Erduran and Keskin \cite{Erduvan2021} solved the two Diophantine equations 

\begin{equation*}
F_k=J_nJ_m  \ \ \textrm{and}
 \  \ J_k=F_nF_m. 
\end{equation*}
In 2023, Gaber \cite{Gaber2023} solved the equations 
\begin{equation*}
P_k=J_n+J_m  \ \ \textrm{and}
 \  Q_k=J_n+J_m , 
\end{equation*}
where $P_k$ and $Q_k$ are the Pell and Pell-Lucas numbers.
\\In this paper, we solved the two Diophantine equations as an extension of those studies

\begin{equation}\label{eqn2}
J_n +J_m =L_k   
\end{equation}

\begin{equation}\label{25}
L_n +L_m =J_k
\end{equation}
in non-negative integers $n,m,\;\text{and}\; k$. More precisely, we proved the following theorems.

\begin{theorem}\label{2}
The non-negative integer solutions $(n,m,k)$, with $n \geq m$, of the Diophantine equation (\ref{eqn2}) are
\[
\left\{
\begin{array}{cccccccc}
(1,1,0), & (2,1,0), & (2,2,0), & (1,0,1), &(2,0,1), &(3,0,2), &(3,1,3),& (3,2,3),

(5,0,5)

\end{array}
\right\}
\]
\end{theorem}

\begin{theorem}\label{26}
The set of non-negative integer solutions $(n,m,k)$, with $n \geq m$, of the Diophantine equation (\ref{25}) are
\[
\left\{
\begin{array}{cccccccc}
(1,0,3), & (2,0,4), & (3,1,4), & (4,3,5), &(6,2,6), & 
(15,1,12)
\end{array}
\right\}
\]
\end{theorem}
Our approach makes use of one of Matveev's most well-known variations \cite{Matveev2000} of Baker's theory to limit each implicit variable to a single primary variable. Matveev's theorem will enable us to obtain the highest values possible for $k$ and $n$. Yet, this value is really high. Then, To reduce these upper boundaries, we use the reduction lemma due to Dujella and Peth\H{o} \cite{Dujella1998}. In the end, we do Sage calculations to identify every answer.\\

This document is arranged as follows: In Section 2, we go over some information that is essential to our proof.
 Section 3 presents the proof of Theorem \ref{2}. Section 4 presents the proof of Theorem \ref{26}.

\section{Preliminary Results}

\subsection{Lucas sequence}
The Binet formula for the Lucas  sequence: 

\begin{equation}
L_n={\theta^n+\eta^n},\; n\geq 0,
\end{equation}
where
\begin{equation*}
\theta=\frac{1+\sqrt{5}}{2},\; \hspace{0.25cm} \;\eta= \frac{1-\sqrt{5}}{2},
\end{equation*}
are the roots of the auxiliary equation $y^2-y-1=0$. The following inequality will assist us in determining the relationships among the three variables $n,m,\;\text{and}\; k$.
\begin{equation}\label{eqn1}
    \theta^{n-1} \leq L_n \leq 2\theta^{n},\; n\geq1.
\end{equation}

\subsection{Jacobsthal sequence }
Jacobsthal sequence has the characteristic equation

\begin{equation*}
y^2-y-2=0,
\end{equation*}
with roots
\begin{equation*}
\theta=2,\; \hspace{0.25cm} \;\eta=-1,
\end{equation*}
and its Binet formula is given by 

\begin{equation}
J_n=\frac{\theta^n-\eta^n}{3}=\frac{2^n-(-1)^n}{3},\; n\geq 0.
\end{equation}
Forthermore, we have 
\begin{equation}\label{eqn4}
2^{n-2} \leq J_n \leq 2^{n-1},\;n\geq 1.
\end{equation}

\subsection{Linear forms in logarithms}

Let $\alpha$ be an algebraic number of degree $d$ with minimal polynomial 
\begin{equation*}
a_0 y^d +a_1 y^{d-1} + \cdots + a_d = a_0 \prod_{i=1}^{d}\left(y-\alpha^{(i)}\right)\;\in\mathbb{Z}[y],
\end{equation*}
where $(a_i , a_j)=1$ for all $a_i \neq a_j$, $a_0>0$ and $\alpha^{(i)}$'s are the conjugates of  $\alpha$. The logarithmic height $h$ of $\alpha$ is defined by 

\begin{equation*}
h(\alpha)= \frac{1}{d} \left(\log {a_0} + \sum_{i=1}^{d}\log\left(max\left\{\left\lvert\alpha^{(i)}\right\lvert,1\right\}\right)\right).
\end{equation*}
The function $h$ has the following properties,
\begin{eqnarray*}
h\left(\alpha_{1}\pm \alpha_{2}\right)&\leq& h\left(\alpha_{1}\right)+h\left(\alpha_{2}\right)+\log 2;\\
h\left(\alpha_{1}\alpha_{2}^{\pm 1}\right)&\leq& h\left(\alpha_{1}\right)+h\left(\alpha_{2}\right);\\
h\left(\alpha^{s}\right)&=&\left\lvert s\right\lvert  \ h\left(\alpha\right) \hspace{0.2cm} \left(s\in \mathbb{Z}\right).
\end{eqnarray*}
These characteristics will be applied in the upcoming sections without particular references.\\
 
We use the following theorem, due to Matveev, in getting upper bounds for the our variables.

\begin{theorem}(Matveev's Theorem) in \cite{Matveev2000}\label{29}
Let L be a real algebraic number field of degree D,  $\alpha_1,\ldots,\alpha_t$ be positive real algebraic numbers in L, $b_1,\ldots,b_t$ be non zero integers and 
\begin{equation*}
\alpha =\alpha_1^{b_1} \cdots \alpha_t^{b_t} -1 \neq 0.
\end{equation*}
Then
\begin{equation}
 \log \left\lvert\alpha\right\lvert > -1.4 \times 30^{t+3} \times t^{4.5} \times D^2 \: (1+\log D) (1+\log B)\: A_1\:A_2\:A_3,
 \end{equation}
where
\begin{equation*}
B\geq \max\{\left\lvert b_1\right\lvert,\ldots,\left\lvert b_t\right\lvert\},    
\end{equation*}
and $A_i \geq \max\{D \: h(\alpha_i) , \left\lvert\log\alpha_i\right\lvert , 0.16\}$ for all $i=1,\ldots,t.$
\end{theorem}

The following Lemma, which is a variation of a lemma by Davenport \cite{Baker1969},  demonstrated by Dujella and Peth\H{o} \cite{Dujella1998}. We will apply this lemma to reduce the upper bounds obtained by Theorem \ref{29}.

\begin{lemma}\label{1}
Let $M$ be a positive integer. Let $\alpha,\mu,A >0, B>1$ be given real numbers. Assume that $\frac{p}{q}$ is a convergent of $\alpha$ such that $q >6M$ and $\epsilon : = \left\|\mu q\right\|-M\left\|\alpha q\right\|>0.$ If $(n,m,\omega) $ is a positive solution to the inequality 
\begin{equation*}
    0<\left\lvert n\alpha -m+\mu\right\rvert<\frac{A}{B^{\omega}},
\end{equation*}
with $n \leq M,$ then 
\begin{equation*}
    \omega < \frac{\log \left(\frac{Aq}{\epsilon}\right)}{\log B}.
\end{equation*}
\end{lemma}
The following Theorem is called \textbf{Legendre Theorem of continued fractions.} More details about this result can be found in \cite{Ddamulira2020}.
\begin{theorem}\label{thm20}
Let $x \in\mathbb{R}$, $ p, q \in \mathbb{Z}$ and let $x=[a_{0},a_{1},...]$.
If
\begin{equation*}
    \left\lvert\frac{p}{q}-x\right\rvert<\frac{1}{2 q^2}
\end{equation*}
then $\frac{p}{q}$ is a convergent of the  continued fraction of $x$. Furthermore, let M and n be nonnegative integers such that $q_n > M$. Put $b=\max\lbrace{a_i:0\leq i\leq n\rbrace}$ then, 
\begin{equation*}
    \frac{1}{\left(b+2\right)s^2}<\left\lvert\frac{r}{s}-x\right\rvert
\end{equation*}
holds for all pairs $(r,s)$ of positive integers with $0<s<M$. 
\end{theorem}

\section{Proof of Theorem 1.1}
\begin{proof}
Assume that $n \geq m$. With Sage, we discovered that the only answers to equation (\ref{eqn2}) in the range $0 \leq m \leq n \leq200$, are the solutions listed above in Theorem \ref{2}. We now assume that $n>200$.

Applying the inequalities (\ref{eqn1}) and (\ref{eqn4}), we obtain
\begin{equation}
    \theta^{k-1} \leq L_k \leq 2^{n} \hspace{0.2cm}  \text{and}\hspace{0.2cm}  2^{n-3} \leq L_k \leq \theta^{k-1}. 
\end{equation}
Then, 
\begin{equation}
    (n-3)\frac{\log 2}{\log \theta} \leq k \leq n \frac{\log 2}{\log \theta}+1.
\end{equation}
Hence,  
\begin{equation}
          n < k < 2n.
\end{equation}

\subsection{Determining the upper boundaries of $k$ and $n$ }
By writing the equation (\ref{eqn2}) as
\begin{equation*}
\frac{2^n}{3} - \theta^k = -J_m + \eta^k + \frac{(-1)^n}{3},
\end{equation*}
we obtain,
\begin{equation}\label{eqn7}
\left\lvert \frac{2^n}{3} - \theta^k\right\lvert <2^m + \frac{11}{6},
\end{equation}
since $\frac{\left\lvert\eta\right\lvert^k}{\sqrt{5}} < \frac{1}{2}$ for all $k\geq0$, and $J_m<2^m$, for all $m\geq0$.
Dividing (\ref{eqn7}) by $2^n/3$, we obtain 
\begin{equation}\label{eqn8}
\left\lvert1-2^{-n} \; \theta^k \; \ 3 \right\lvert < \frac{9}{2^{n-m}}.
\end{equation}
Thus,  
\begin{equation}\label{eqn1111}
    \log\left\lvert1-2^{-n} \; \theta^k \; \ 3\right\lvert < \log 9 -(n-m)\log2.
\end{equation}
\\Define
$$\zeta_1:=2^{-n} \; \theta^k \; \ 3 -1 .$$
If $\zeta_1= 0$, then $3 \theta^k=  2^n$.
By conjugating the last equation, we obtain $3 \eta^k=2^n$.Taking the absolute values, we obtain that $2^n=\left\lvert3 \eta^k\right\lvert<5$ which is a contradiction. Hence $\zeta_1 \neq 0$.\\

Assume that
$$\alpha_1=2,\; \alpha_2=\theta,\; \alpha_3=3,\  b_1=-n,\; b_2=k,\; b_3=1, \textrm{and} \ t=3.$$
Since $\max\{\left\lvert-n\right\lvert,\left\lvert k\right\lvert,\left\lvert 1\right\lvert\}=k$, we consider $B=k$. The smallest field contains $\alpha_1,\;\alpha_2,\; \text{and} \; \alpha_3$ is $\mathbb{Q}$$(\sqrt{5})$. Consequently, $D=2$. The logarithmic heights are 
\begin{equation*}
 h(\alpha_1)=\log 2,   
\end{equation*}
\begin{equation*}
h(\alpha_2)=(\log \theta)/2, 
\end{equation*}
and
\begin{equation*}
 h(\alpha_3)=\log 3.
\end{equation*}
So, we take 
\begin{equation*}
A_1=1.4,\; A_2=0.5,\; \text{and}\; A_3=2.2.
\end{equation*}
Applying Theorem \ref{29}, we obtain 
\begin{equation}\label{eqn9}
    \log\left\lvert\zeta_1\right\lvert > -C(1+\log k),
\end{equation}
where $C=1.54 \times 1.4 \times 30^6 \times 3^{4.5} \times 2^2\times (1+\log2) < 14.938\times 10^{11}$.\\

Using (\ref{eqn1111}) and  (\ref{eqn9}), we get
\begin{equation*}
\log 9 -(n-m)\log2 > -14.938 \times 10^{11} \times (1 + \log k) > -2.99 \times 10^{12} \times \log k, 
\end{equation*}
hence 
\begin{equation}\label{eqn13}
(n-m) \log2 < 3 \times 10^{12}\times \log k. 
\end{equation}
We now write the equation (\ref{eqn2}) as
\begin{equation*}
\frac{2^n}{3}+\frac{2^m}{3}-\theta^k=\eta^k+\frac{(-1)^n}{3}+\frac{(-1)^m}{3}.
\end{equation*}
Consequently,
\begin{equation*}
\left\lvert\frac{2^n}{3} \left(1+2^{m-n}\right) - \theta^k\right\lvert < \frac{11}{6}.
\end{equation*}
By dividing the above inequality by $\frac{2^n}{3} \left(1+2^{m-n}\right)$ , we obtain
\begin{equation}\label{eqn10}
\left\lvert1-\theta^k\; 2^{-n}\; \  3 \left(1+2^{m-n}\right)^{-1}\right\lvert < \frac{11}{2}\times \frac{1}{2^n}.  
\end{equation}
Thus,
\begin{equation} \label{eqn1010}
\log \left\lvert1-\theta^k\; 2^{-n}\; \  3 \left(1+2^{m-n}\right)^{-1}\right\lvert < \log\frac{11}{2} - n\log2. 
\end{equation}
Define
$$\zeta_2:=\theta^k\; 2^{-n}\; \  3 \left(1+2^{m-n}\right)^{-1} -1 .$$
If $\zeta_2=0$, then $3 \theta^k=\left(2^n+2^m\right)$. By conjugating the last equation, we obtain that $\left(2^n+2^m\right)=\left\lvert3 \eta^k\right\lvert<5$, which is a contradiction. Hence $\zeta_2 \neq 0$.\\

Assume that
$$\alpha_1=2,\ \alpha_2=\theta,\ \alpha_3=  3
\left(1+2^{m-n}\right)^{-1},\ b_1=-n,\ b_2=k,\ b_3=1, \ \text{and} \ t=3.$$
We can see that $\alpha_1$,\;$\alpha_2$, and $\alpha_3$ are in $\mathbb{Q}$$(\sqrt{5})$, so we take $D=2$. Also, we take $B=k$. The logarithmic heights are 
\begin{equation*}
 h(\alpha_1)=\log 2,   
\end{equation*}
\begin{equation*}
h(\alpha_2)=\frac{1}{2} \log \theta, 
\end{equation*}
and
\begin{eqnarray*}
 h (\alpha_3)=h\left(\  3 \left(1+2^{m-n}\right)\right)^{-1})&\leq &h\left(\  3\right) + h\left(2^{m-n}\right) + \log2 \\ &\leq &\log3 + \left\lvert m-n\right\lvert\: h(2)+\log2 
\\&=&\log6 + (n-m) \log2.
\end{eqnarray*} 
Since
\begin{equation*}
\left\lvert\log\alpha_3\right\lvert = \left\lvert\log  \  3 - \log(1+2^{m-n})\right\lvert<\log \  3+\log2< 1.8,
\end{equation*}
therefore, we take 
\begin{equation*}
A_1=1.4,\; A_2=0.5,\; \text{and}\;A_3=4 + 2(n-m) \log 2\;.
\end{equation*}
By applying Theorem \ref{29}, we obtain  
\begin{equation*}
\log \left\lvert\zeta_2\right\lvert > -C (1+\log k) \times (4 + 2(n-m) \log2), 
\end{equation*}
where $C=(1.4)^{2} \times .5 \times 30^6 \times 3^{4.5} \times 2^2\times (1+\log2)< 6.79\times10^{11}$.\\
Therefore,
\begin{equation}\label{eqn11}
\log \left\lvert\zeta_2\right\lvert > -6.79 \times 10^{11} \times (1+ \log k) > -1.36 \times 10^{12} \times (\log k)\times (4 + 2(n-m) \log2).
\end{equation}
Using (\ref{eqn1010}) and (\ref{eqn11}), we obtain
\begin{equation}\label{eqn12}
n\log2 <1.37 \times 10^{12}\times (\log k) \times(4 + 2(n-m) \log2).
\end{equation}
Substituting  by inequality (\ref{eqn13}) in inequality (\ref{eqn12}), we obtain 
\begin{equation}\label{eqn14}
n \log2 < 1.37 \times10^{12}\times (\log k) \left(4 +6 \times 10^{12} \log k \right). 
\end{equation}
Using (3.2) and (\ref{eqn14}) , we get
\begin{equation}\label{eqn37}
k \log2 < 2.74 \times10^{12}\times (\log k) \left(4 +6 \times 10^{12} \log k\right). 
\end{equation}
Solving inequality (\ref{eqn37}) by  Sage, we obtain 
\begin{equation}\label{40}
n< k <1.1 \times 10^{29}. 
\end{equation}
\subsection{Lowering the limit on $n$ }
Currently, Lemma \ref{1} will be applied to lower the upper bound in inequality (\ref{40}) .\\
Let
\begin{equation}\label{eqn17}
Q = k \log\theta - n \log2   + \log 3.
\end{equation}
Rewrite equation (\ref{eqn8}) as
\begin{equation}
\left\lvert e^Q - 1\right\lvert < \frac{9}{2^{n-m}}.
\end{equation}
As we proved that $\zeta_1 \neq 0$, we obtain that $Q \neq 0$. If $Q > 0$, we obtain that
\begin{equation}
0 < Q < e^Q - 1 =\left\lvert e^Q - 1\right\lvert < \frac{9}{2^{n-m}}.
\end{equation}
If $Q < 0$, and $n-m \geq 25$, then
\begin{equation*}
\left\lvert e^Q -1\right\lvert=1-e^Q<\frac{1}{2},
\end{equation*}
therefore
\begin{equation*}
e^{\left\lvert Q\right\lvert}=e^{-Q} < {2}.
\end{equation*}
Therefore, we obtain 
\begin{equation}\label{eqn18}
0 < \left\lvert Q\right\lvert < e^{|Q|} - 1 = e^{\left\lvert Q\right\lvert} (1-e^Q) = e^{\left\lvert Q\right\lvert} \left\lvert e^Q -1\right\lvert< \frac{18}{\  2^{n-m}}.
\end{equation}
In both cases $(Q>0\; \text{and}\; Q<0)$, inequality (\ref{eqn18}) is true. By substituting the formula (\ref{eqn17}) for $Q$ in the inequality above and dividing the result by $\log 2$, we obtain
\begin{equation}\label{eqn19}
0 < \left\lvert k\frac{\log \theta}{\log2} - n + \frac{\log\  3}{\log2}\right\lvert < \frac{26}{2^{n-m}}.
\end{equation}
The next thing we use Lemma \ref{1}  on inequality (\ref{eqn19}). Let
\begin{equation*}
\mu = \frac{\log\ 3}{\log 2},\hspace{1cm} \alpha = \frac{\log\theta}{\log2},\hspace{1cm} w = n-m,\hspace{1cm} A = 26,\hspace{0.5cm}\text{and}\hspace{0.5cm} B=2 .
\end{equation*}
 Using inequality (\ref{40}), we are able to take $M =3\times10^{29}$ as an upper bound on $k$. For 
$$q_{75} = 252339790309653189029774211371593442> 6M,$$
where the denominator of the $75^{th}$ convergent of the continued fraction of $\alpha$ is $q_{75},$ we found that 
$$\epsilon = \left\|\mu \:q_{75} \right\|-M \left\|\alpha \: q_{75}\right\|>0.140378035627 .$$
 Lemma \ref{1} implies
\begin{equation}
n-m < 127.
\end{equation}

To handle the previous upper bound of $n$ in equality (\ref{40}), we will get to work on inequality (\ref{eqn10}).\\ Let
\begin{equation}\label{eqn22}
S = k \log\theta - n \log2 + \log  3 \left(1+2^{m-n}\right)^{-1}. 
\end{equation}
Inequality (3.8) implies that 
\begin{equation}
\left\lvert 1-e^S\right\lvert< \frac{11}{2^{n+1}}.
\end{equation}
Observe that $S \neq 0$, since $\zeta_2\neq 0$. If $S > 0$, we obtain that  
\begin{equation}
0 < S < e^S - 1 =\left\lvert e^S - 1\right\lvert < \frac{11}{2^{n+1}},
\end{equation}
and if $S < 0$, we have
\begin{equation}
1-e^S=\left\lvert e^S-1\right\lvert  < \frac{11}{2^{n+1}}< \frac{1}{2} \;\: \; \text{for all}\; n > 200,
\end{equation}
the last inequality implies that\;$e^{\left\lvert S\right\lvert} <{2} $\;. Hence
\begin{equation}\label{eqn21}
0 < \left\lvert S\right\lvert < e^{\left\lvert S\right\lvert} - 1 = e^{\left\lvert S\right\lvert} \left(1-e^S\right) = 
e^{\left\lvert S\right\lvert} \left\lvert e^S -1\right\lvert < \frac{11}{2^n}.
\end{equation}
In both cases, ($S>0$, or $S<0$), inequality (\ref{eqn21}) is true. By substituting the formula (\ref{eqn22}) for $S$ in the inequality above and dividing the result by $\log 2$, we obtain
\begin{equation}
0 < \left\lvert k\frac{\log \theta}{\log2} - n + \frac{\log\  {3}{\left(1+2^{m-n}\right)^{-1}}  \ }{\log2}\right\lvert < \frac{1}{2^{n-4}}.
\end{equation}
Now, we apply Lemma \ref{1} on the last inequality. Let
\begin{equation*}
 \mu =\frac{\log{3} \left(1+2^{-\left(n-m\right)}\right)^{-1}}{\log2},\hspace{0.5cm} \alpha = \frac{\log \theta}{\log2},\hspace{0.5cm} w = n-4,\hspace{0.5cm}A = 1,\hspace{0.5cm}\text{and}\hspace{0.5cm} B = 2.
\end{equation*}
Using inequality (\ref{40}), we can take $M=3 \times 10^{29}$ as an upper bound on $k.$ For
$$q_{75}=252339790309653189029774211371593442 > 6M,$$ where the denominator of the $75^{th}$ convergent of the continued fraction of $\alpha$ is $q_{75},$ we found that,  for  all integers $n-m \in[0,128],$ the smallest $\epsilon$ occurs at $n-m= 121   $  except if $n - m = 1$
$$\epsilon = \left\|\mu \:q_{75} \right\|-M \left\|\alpha \: q_{75}\right\| > 0.00343788493     .$$
By applying Lemma \ref{1}, we obtain
\begin{equation}\label{eqn23}
n < 200,
\end{equation}
 This contradicts the assumption, which is $n > 200$.  therefore, equation (\ref{eqn2}) has no solution when $n > 200.$
 Thus, we obtain  $n\leq 200 $ so
 \begin{equation}
   m\leq n \leq 200 
 \end{equation}

 and
 \begin{equation}
 k<2n<400.
 \end{equation}
If $n-m=1$ we obtain $\epsilon$ always negative. So if ${n-m}=1 $. Equation (\ref{eqn2}) becomes $L_k =J_{n}+J_{n+1}$, which can be reduced to
 \begin{equation}\label{eqn62}
     L_{k}=2^{n}.
 \end{equation}
 Then $k<2n$ and  from equation (\ref{40}) we obtain $n<3\times 10^{29}$. Using the same method of verification $\zeta _{1}<\frac{9}{2^{m-n}}$, we obtain
 \begin{equation*}
     \theta^{k}2^{-n}-1<\frac{3}{2^{n}},
 \end{equation*}
 thus
 \begin{equation*}
     \left\lvert k\frac{\log \theta}{\log 2}-n\right\rvert<\frac{18}{2^n}.
 \end{equation*}
$\frac{18}{2^n}<\frac{1}{2k^2}$ for all $n > 10$ since $k<2n$.
 From Theorem \ref{thm20}, this implies that $\frac{n}{k}$ is a convergent of $\frac{\log \theta}{\log 2}$. We know that $k<M$ then using  computation with sage gives that 
 \begin{equation*}
     q_{67}\leq M < q_{68} \hspace{0.3cm} \text{and} \hspace{0.3cm} b:=\max \{a_{i}: i=0,2,...,69\}=134.
 \end{equation*}
Therefore, Theorem \ref{thm20} implies that
 \begin{equation*}
     \frac{1}{\left(b+2\right) k}<\left\lvert k\frac{\log \theta}{\log 2}-n\right\rvert< \frac{18}{2^n}.
 \end{equation*}
 Thus 
 \begin{equation*}
     2^{n}<1.944 \times 10^{32}.
 \end{equation*}
 Then $n\leq 110.$ Hence the Theorem \ref{eqn2} is proved.
\end{proof}
\section{Proof of Theorem 1.2}
\begin{proof}
By using Sage, It was discovered that the aforementioned solutions in Theorem \ref{26} are the only solutions of equation (\ref{25}) inside of the range $0 \leq m\leq
n \leq 200.$  As a result, suppose that $n>200$.

\subsection{Finding relationship between $k$ and $n$.}
Applying the inequalities (\ref{eqn1}) and (\ref{eqn4}) to create a relation between $k$ and $n$, we obtain
\begin{equation}
    \ 2^{k-2} \leq J_k \leq 2L_n\leq 4\theta^{n},
\end{equation}  
then,
\begin{equation}
    \  k \leq n \frac{\log\theta}{\log 2}+4.
\end{equation}
Hence,
\begin{equation}
             k < n.
\end{equation}

\subsection{Finding upper bound of $n$}
By writing equation (\ref{25}) as 
\begin{equation*} 
\theta^{n}-\frac{2^k}{3}=-L_m-\frac{(-1)^k}{3}-\eta^n,
\end{equation*}
we obtain,
\begin{equation}\label{eqn55}
\left\lvert \theta^{n}-\frac{2^k}{3}\right\lvert\leq \left\lvert \ L_m \right\lvert\ +{\left\lvert\eta\right\lvert^n}+\frac{1}{3}<2\theta^m+\frac{11}{6},
\end{equation}
since $\frac{\left\lvert\eta\right\lvert^n}{\sqrt{5}}<\frac{1}{2}$  for all $n\geq0$, and $L_m<2\theta^m$  for every $m\geq 1$. Dividing (\ref{eqn55}) by $\theta^{n}$, we obtain
\begin{equation}\label{28}
\left\lvert1-2^k\; \theta^{-n}\; \frac{1}{3}\right\lvert<\frac{4}{\theta^{n-m}}.
\end{equation}
As a result,
\begin{equation}\label{200}
\log \left\lvert1-2^k\; \theta^{-n}\; \frac{1}{3}\right\lvert<\log4 - \log\theta^{n-m}.
\end{equation}
Define
$$\zeta_3:=2^k\; \theta^{-n}\; \frac{1}{3} - 1.$$
If $\zeta_3 = 0$, we obtain that $3\theta^n=2^K$. By conjugating the last equation, we obtain that $2^k =3\left\lvert\eta^n\right\lvert<5$, which is a contradiction. Hence $\zeta_3 \neq 0$.\\

Let
$$\alpha_1=2,\ \alpha_2=\theta,\ \alpha_3=\frac{1}{3}, \ b_1=k,\ b_2=-n,\ b_3=1,\;\text{and}\;t=3.$$
Since $n>k$, we consider $B=n$, and  $L=\mathbb{Q}(\sqrt{5})$, which has a degree $D=2$. The logarithmic heights are  
\begin{equation*}
 h(\alpha_1)=\log 2,   
\end{equation*}
\begin{equation*}
h(\alpha_2)=\frac{\log \theta}{2},                
\end{equation*}
and
\begin{equation*}
 h(\alpha_3)=h \left(\frac{1}{3}\right)=\log {3}
\end{equation*}
So, we take
\begin{equation*}
A_1=1.4,\; A_2=0.5,\; \textrm{and}\; A_3=2.2\;.
\end{equation*}
Applying Theorem \ref{29}, we obtain 
\begin{equation*}
\log\left\lvert\zeta_3\right\lvert>-C (1+\log n),
\end{equation*}
where $C=1.54 \times 1.4 \times 30^6 \times 3^{4.5} \times 2^2\times (1+\log2) < 14.938 \times 10^{11}$.
Using (\ref{200}) and (\ref{29})  we obtain
\begin{equation*}
\log4 - \log\theta^{n-m} > -14.938 \times 10^{11}(1+ \log n) > -2.987 \times10^{12} \times \log n,
\end{equation*}
thus,
\begin{equation}\label{111}
(n-m )\log\theta < 3 \times10^{12} \log n .
\end{equation}
We now write equation (\ref{25}) as
\begin{equation*}
\theta^n + \theta^m -\frac{2^k}{3}= -\eta^n - \eta^m-\frac{(-1)^k}{3}
.\end{equation*} 
Thus,
 \begin{equation}\label{equatiom1}
      \left\lvert\theta^n \left(1+\theta^{m-n}\right) - \frac{2^k}{3}\right\lvert \leq \left\lvert\eta|^{n}+|\eta\right\lvert^{m}+\frac{1}{3}<4,
\end{equation}
where $\frac{\left\lvert\eta\right\lvert^n}{\sqrt{5}}<\frac{1}{2}$ for all $n\geq 0$. Dividing (\ref{equatiom1}) by  $\theta^n \left(1+\theta^{m-n}\right)$, we obtain 
\begin{equation}\label{30}
\left\lvert1-2^k \;\theta^{-n}\; \frac{1}{3} \left(1+\theta^{m-n}\right)^{-1}\right\lvert < \frac{4}{\theta^n}.
\end{equation}
Thus,
\begin{equation}\label{201}
\log \left\lvert1-2^k \;\theta^{-n}\; \frac{1}{3} \left(1+\theta^{m-n}\right)^{-1}\right\lvert <\log4 - n \log\theta.
\end{equation}
Define

$$\zeta_4:=2^k \;\theta^{-n}\; \frac{1}{3} \left(1+\theta^{m-n}\right)^{-1} - 1.$$
If $\zeta_4=0$, then $3\left(\theta^n + \theta^m\right) = \ 2^k$. By conjugating the last equation, we obtain $2^k =3\left\lvert\eta^n + \eta^m\right\lvert<9$, which is a contradiction. Hence Hence $\zeta_4 \neq 0$. \\

Let
$$\alpha_1=2,\; \alpha_2=\theta,\; \alpha_3=\frac{1}{3} \left(1+\theta^{m-n}\right)^{-1}, \ b_1=k,\; b_2=-n,\; b_3=1,\;\text{and}\;t=3.$$
We take $B=n$, $L=\mathbb{Q}(\sqrt{5})$, and $D=2$. The logarithmic heights are 
\begin{equation*}
 h(\alpha_1)=\log 2,   
\end{equation*}
\begin{equation*}
h(\alpha_2)=\frac{\log \theta}{2}, 
\end{equation*}
\begin{eqnarray*}
h\left(\frac{1}{3} \left(1+\theta^{m-n}\right)^{-1}\right)&\leq& h\left(\frac{1}{3}\right) + h\left(1+\theta^{m-n}\right) \\&\leq& \log3 + h\left(\theta^{m-n}\right) + \log2 \\&=& \log6 + \left\lvert m-n\right\lvert \;h\left(\theta\right)\\ &=& \log6 + (n-m)\; \frac{\log \theta}{2}.
\end{eqnarray*}
Since
\begin{eqnarray*}
\left\lvert\log\alpha_3\right\lvert&=&\left\lvert\log\frac{1}{3}-\log \left(1+\theta^{m-n}\right)\right\lvert\\&<&\left\lvert\log\frac{1}{3}\right\lvert + \left\lvert\log\left(1+\theta^{m-n}\right)\right\lvert 
\\&<&\left\lvert\log\frac{1}{3}\right\lvert + \left\lvert\log2\right\lvert < 1.1 + 0.7=1.8,
\end{eqnarray*} 
hence, we can take
\begin{equation*}
A_1=1.4,\; A_2=0.5,\;\text{and}\; A_3=4 + (n-m) \log\theta >\max \left\{2h(\alpha_3) , \left\lvert\log\alpha_3\right\lvert , 0.16 \right\}.
\end{equation*}
Applying Theorem \ref{26}, we obtain
\begin{equation}\label{202}
\log \left\lvert\zeta_4\right\lvert > - C\: (1+\log n) \times((n-m) \log\theta +4), 
\end{equation}
 where $C=(1.4)^{2} \times 0.5 \times 30^6 \times 3^{4.5} \times 2^2 \times (1+\log2) < 6.79 \times 10^{11}$. Using (\ref{201}) and (\ref{202}), we obtain 
\begin{equation*}
\log4 - n \log\theta >-6.79 \times 10^{11} \times (1 + \log n) > -1.36 \times 10^{12} \times \log n \times \left(\left(n-m\right) \log\theta +4 \right),
\end{equation*}
which implies that 
\begin{equation}\label{34}
n \log \theta < 1.37 \times 10^{12} \times \log n \times \left((n-m) \log\theta +4\right).
\end{equation}
Using (\ref{111}) and (\ref{34}), we obtain   
\begin{equation}\label{equatiom2}
n \log \theta < 1.37 \times 10^{12} \times \log n \times (4+3\times10^{12} \times \log n).
\end{equation}
When Sage solves the inequality (\ref{equatiom2}), we obtain  
\begin{equation}\label{31}
n<3.7 \times10^{28}.
\end{equation}
\subsection{Reducing the bound on $n$}
Currently, we will lower the obtained upper bound in inequality (\ref{31}) using Lemma \ref{1}. Let
\begin{equation}\label{33}
 Q = k \log2 - n \log\theta + \log\frac{1}{3}.     
\end{equation}
Inequality (\ref{30}) implies that
\begin{equation*}
\left\lvert1-e^Q\right\lvert<\frac{4}{\theta^{n-m}}.
\end{equation*}
Note that $Q \neq 0$, since $\zeta_3\neq 0$. If $Q > 0$, we obtain 
\begin{equation}
0 < Q < e^Q - 1 =\left\lvert e^Q - 1\right\lvert < \frac{4}{\theta^{n-m}}.
\end{equation}
If $Q < 0$, we suppose that\; $n-m \geq 20$, and since  $\frac{4}{\theta^{n-m}} < \frac{1}{2}$ holds for  $n-m \geq 20$, we get

\begin{equation*}
\left\lvert1-e^Q\right\lvert =1-e^Q < \frac{1}{2}.
\end{equation*}
Consequently, 
\begin{equation*}
e^{\left\lvert Q\right\lvert} < 2.
\end{equation*}
Therefore,
\begin{equation}\label{32} 
0 < \left\lvert Q\right\lvert < e^{\left\lvert Q\right\lvert} - 1 = e^{\left\lvert Q\right\lvert} (1-e^Q) = e^{\left\lvert Q\right\lvert}\; \left\lvert e^Q-1\right\lvert < \frac{8}{\theta^{n-m}}. 
\end{equation}
In both cases $(Q>0 \ \text{and} \ Q<0)$, inequality (\ref{32}) is true. By substituting the formula (\ref{33}) for $Q$ in the inequality above and dividing the result by $\log 2$, we obtain
\begin{equation}
0 < \left\lvert n\frac{\log \theta}{ \log \ 2} - k + \frac{\log \ 3}{\log \ 2}\right\lvert < \frac{12}{\theta^{n-m}}. 
\end{equation}

Let

\begin{equation*}
 \mu = \frac{\log\ {3}}{\log 2},\hspace{0.5cm} \alpha = \frac{\log \theta}{\log 2},\hspace{0.5cm} w = n-m,\hspace{0.25cm} A = 12,\hspace{0.5cm} \text{and} \hspace{0.25cm}  B = \theta .
\end{equation*}
Using inequality (\ref{31}), we can take $M=4$×$10^{28}$ as an upper bound of $n.$ For 
$$q_{70} = 228666343422267608843910896109913 > 6M $$

where the denominator of the $70^{th}$ convergent of the continued fraction of $\alpha$ is $q_{70},$ we found that $$\epsilon = \left\|\mu \:q_{70} \right\|-M \left\|\alpha \: q_{70}\right\| >0.0328403974748. $$ By applying Lemma \ref{1}, we obtain 
\begin{equation}\label{205}
 n-m <168.
\end{equation}

Now we consider 
\begin{equation*}
\left\lvert1-2^k \theta^{-n} \frac{1}{3} \left(1+\theta^{m-n}\right)^{-1}\right\lvert < \frac{4}{\theta^n}.
\end{equation*}
Let
\begin{equation}\label{35}
S = k \log2 - n \log\theta + \log \frac{1}{3} \left(1+\theta^{m-n}\right)^{-1}.
\end{equation}
Then
\begin{equation} 
\left\lvert1-e^S\right\lvert < \frac{4}{\theta^n}. 
\end{equation}
Note that $S \neq 0$ since $\zeta_4\neq0$. If $S > 0$, we obtain 
\begin{equation}
0 < S < e^S - 1 = \left\lvert e^S - 1\right\lvert < \frac{2}{\theta^n}.
\end{equation}
Since $n > 200$, we have $\frac{2}{\theta^n} < \frac{1}{2}$, which implies that for $S<0$, we have  
\begin{equation*}
1- e^S=\left\lvert e^S -1\right\lvert < \frac{1}{2}.
\end{equation*}
It follows that
\begin{equation*}
e^{\left\lvert S\right\lvert}=e^{-S} < 2.
\end{equation*}
Thus,
\begin{equation}\label{112}
0 < \left\lvert S\right\lvert < e^{\left\lvert S\right\lvert} - 1 = e^{\left\lvert S\right\lvert} \left(1-e^S\right) = e^{\left\lvert S\right\lvert} \left\lvert e^S -1\right\lvert < \frac{8}{\theta^n}.
\end{equation}
In both cases  ($S>0,S<0$), inequality (\ref{112}) is true. By substituting the formula (\ref{35}) for $S$ in the inequality above and dividing the result by $\log 2$, we obtain  
\begin{equation}
0 < \left\lvert k - n \frac{\log\theta}{\log 2} + \frac{\log \frac{1}{3} (1+\theta^{m-n})^{-1}}{\log 2} \right\lvert< \frac{12}{\theta^{n}}.
\end{equation}
\begin{equation}
0 < \left\lvert   n \frac{\log\theta}{\log 2} - k  + \frac{\log 3 (1+\theta^{m-n})}{\log 2} \right\lvert< \frac{12}{\theta^{n}}.
\end{equation}

Let
\begin{equation*}
 \mu = \frac{\log 3 (1+\theta^{m-n})}{\log 2},\hspace{0.5cm} \alpha = \frac{\log\theta}{\log 2},\hspace{0.25cm}\hspace{0.5cm} w = n,\hspace{0.25cm} A = 12,\hspace{0.5cm} \text{and} \hspace{0.25cm}  B = \theta .
\end{equation*}
Using inequality (\ref{31}), we can take $ M = 4 \times 10^{28} > n > k $ as an upper bound of $k.$ Note that for the same 
$$q_{70} = 228666343422267608843910896109913 > 6M,$$ 
where the denominator of the $70^{th}$ convergent of the continued fraction of $\alpha$ is $q_{70},$  we found that  for  all integers $n-m \in[0,168]$ the smallest $\epsilon$ is 
$$\epsilon = \left\|\mu \:q_{70} \right\|-M \left\|\alpha \: q_{70}\right\| >0..004.$$
 By applying Lemma \ref{1}, we obtain 
\begin{equation}\label{1093}
n <172,
\end{equation}
This contradicts the assumption, which is $n > 200$.
Therefore, there is no solution to equation (\ref{25}) when $n > 200.$ 

\end{proof}

\makeatletter

\vspace{-0.3cm}
\makeatother



\end{document}